\documentclass{gtart}
\gtart

\usepackage{amsthm}
\usepackage{amsgen}
\usepackage{amsmath,amssymb}
\usepackage{curvesls}
\usepackage{epic}
\usepackage[all]{xy}
\usepackage[dvips]{graphicx}
\usepackage{psfrag}

\newtheorem{thm}{Theorem}[section]

\newtheorem{prop}[thm]{Proposition}

\theoremstyle{definition}

\newcommand{\Real}{{\mathbb R }}

\newcommand{\SO}{{\mathrm{SO}}}
\newcommand{\EK}[1]{{\mathrm{EC}({#1})}}

\newcommand{\coba}{{\mathcal{CB}}}
\newcommand{\frdi}{{\mathfrak{f}\mathcal{D}}}

\newcommand{\CD}{{\mathrm{CDiff}}}
\newcommand{\Trans}{{\mathrm{Trans}}}
\newcommand{\Scale}{{\mathrm{Scale}}}

\begin{document}

\title{The framed discs operad is cyclic}

\author{Ryan Budney}

\address{IH\'ES, Le Bois-Marie, 35, route de Chartres \\
F-91440 Bures-sur-Yvette, FRANCE
}

\email{budney@ihes.fr}

\begin{abstract}
The operad of framed little discs is shown to be equivalent to a cyclic 
operad. This answers a conjecture of Salvatore in the affirmative, posed 
at the workshop `Knots and Operads in Roma,' at Universit\`a di Roma 
``La Sapienza'' in July of 2006. 
\end{abstract}


\maketitle

\section{Statement and proof}\label{StatnProof}

Both the author \cite{Bud} and Paolo Salvatore \cite{Sal} have constructed actions of the 
operad of little $2$-cubes on the `framed' long knot spaces $\EK{1,D^n}$ defined in \cite{Bud}
which extend the homotopy-associative product given by the connect-sum operation.  
Computations of the homotopy-type \cite{KS} and homology \cite{Cohen} of $\EK{1,D^2}$ suggest in several
ways (some more vague than others) that the spaces $\EK{1,D^n}$ may admit an action of the operad of 
framed little $2$-discs.  A concrete example: provided such an action exists, it would allow for a more 
compact and `operadic' description of the homotopy-type of $\EK{1,D^2}$ than the one given in \cite{Bud, KS}. 
The author has made several attempts to construct an action of the operad of framed $2$-discs $\frdi_2$ on $\EK{1,D^n}$ in the spirit of `Little cubes and long knots' \cite{Bud}, as of yet without success.  
Salvatore constructs the spaces $\EK{1,D^n}$ for $n \geq 3$ as the totalisation of a multiplicative operad, further developing the work of Sinha \cite{Sin}. The operad Salvatore uses is equivalent to the operad of 
framed little $(n+1)$-discs $\frdi_{n+1}$. Theorems of McClure and Smith imply that the totalisation of a multiplicative operad is a module over the operad of $2$-cubes \cite{MS}. Motivated by these constructions, Salvatore asked if the operad of framed discs is cyclic. The geometric techniques used in the author's 
failed attempts to construct an action of $\frdi_2$ directly on $\EK{1,D^n}$ turn out to answer his
question. In Proposition \ref{CBprop}, $\frdi_n$ is shown to be equivalent in the sense 
of Fiedorowicz \cite{Fied} 
to the operad of  `conformal $n$-balls' $\coba_n$. In Theorem \ref{mainthm}, $\coba_n$ is shown to be cyclic.

This paper follows the conventions of the book `Operads in Algebra, Topology
and Physics' by Markl, Shnider and Stasheff \cite{MSS}. Let $S^n$ be the 
unit sphere in $\Real^{n+1}$.  
Consider $\Real^0 \subset \Real^1 \subset \Real^2 \subset \cdots $. Let $p_+ = (0,\cdots,0,1) \in S^n$ and 
$p_- = -p_+$ be the north and south poles of $S^n$ respectively. 
Let $D^n \subset S^n$ be the `southern' hemisphere in $S^n$, meaning $\partial D^n = S^{n-1}$ and $p_- \in D^n$. 
Given a point $v \in S^n$, identify the tangent space $T_{-v}S^n$ with $S^n \setminus \{v\}$ via 
stereographic projection from $v$. 
Since stereographic projection is a conformal diffeomorphism, one can consider a conformal 
affine-linear transformation of $T_{-v}S^n$ to be a conformal diffeomorphism of $S^n \setminus \{v\}$
and further extend it to a conformal diffeomorphism of $S^n$.  
Let $\CD(S^n)$ be the group of orientation-preserving conformal diffeomorphisms of $S^n$. 
Let $\CD(S^n,v) \subset \CD(S^n)$ be the subgroup corresponding to the conformal affine-linear transformations of $T_{-v}S^n$. 
Let $\Trans(v) \subset \CD(S^n,v)$ be the subgroup consisting only of translations
in $T_{-v}S^n$. 
Let $\SO_{n+1}$ denote the group of orientation-preserving isometries of 
$S^n$. 
Let $\SO_n(v)$ be defined as $\SO_{n+1} \cap \CD(S^n,v)$. 
Define $\Scale(v) \subset \CD(S^n,v)$ to be the scalar multiples of the identity on
$T_{-v}S^n$.  

For $n \geq 2$, $\CD(S^n)$ is known to be a $n+1 \choose 2$-dimensional Lie group,
isomorphic to the group of hyperbolic isometries of hyperbolic $(n+1)$-space.
The isomorphism is given by restriction to the projectivised light-cone of
the Minkowski model.  Stated another way, for any $v \in S^n$, 
$\CD(S^n)$ is generated by the subgroups $\SO_{n+1}$, $\Trans(v)$ and
$\Scale(v)$. Moreover, the stabiliser of $v \in S^n$ under the action of $\CD(S^n)$ 
is $\CD(S^n,v)$. For details, see Schoen and Yau's book 
\cite{SchYau} \S VI Theorem 1.1.  

Let $\pi : S^n \to S^n$ denote rotation by 180-degrees around $S^{n-2} \subset S^n$. 
A $(j+1)$-tuple $(\pi,f_1,\cdots,f_j) \in \CD(S^n)^{j+1}$ will be called $j$ little 
conformal $n$-balls if $f_i(D^n) \subset D^n$ for all $i$, and if the interior of $f_i(D^n)$ is 
disjoint from $f_k(D^n)$ for all $i \neq k$. The set of all such $(j+1)$-tuples 
is denoted $\coba_n(j)$. There are maps 
$\coba_n(i) \times \left( \coba_n(j_1) \times \cdots \coba_n(j_i) \right)
\to \coba_n(j_1 + \cdots + j_i)$ defined by sending
$\left( (\pi,f_1,\cdots,f_i),
  (\pi,g_{1 1},\cdots,g_{1 j_1}), \cdots,
  (\pi,g_{i 1},\cdots,g_{i j_i}) \right)$ to
$\left( \pi, f_1 g_{1 1},\cdots,f_1 g_{1 j_1},\right.$
	$\left.\cdots, f_i g_{i 1} \cdots f_i g_{i j_i} \right)$.
Consider $\coba_n(j)$ to be a subspace of $\CD(S^n)^{\{0,1,\cdots,j\}}$.
Then there is a natural right action of $\Sigma_j$ on $\coba_n(j)$. The above maps 
together with the actions of $\Sigma_j$ for $j \in \{1,2,\cdots\}$ endow $\coba_n$ with 
the structure of an operad, hereby named the operad of conformal $n$-balls.  Observe that 
for $n \geq 2$ the operad of framed little $n$-discs $\frdi_n$ is the suboperad 
of $\coba_n$ such that the north pole is fixed by each disc: 
$f_i(p_+)=p_+$ for all $i \in \{1,\cdots,j\}$.

\begin{prop} \label{CBprop}
There is a deformation-retraction of $\coba_n$ to the suboperad
$\frdi_n$. Moreover, the inclusions 
$\frdi_n(j) \to \coba_n(j)$ are $\Sigma_j$-equivariant
homotopy-equivalences for all $j \in {\mathbb N}$.
\begin{proof}
To construct the homotopy-inverse $\coba_n \to \frdi_n$,
consider $f \in \CD(S^n)$ such that $f(D^n) \subset D^n$.
$f$ can be written uniquely as a composite $f =  Y \circ S \circ M \circ Q$
where $Y \in \Trans(p_+)$, $S \in \Scale(p_+)$, $M \in \SO_n(p_+)$, $Q \in \Trans(p_-)$
and $f_i(p_-) = Q(p_-)$:
\begin{itemize}
\item $Q^{-1}$ is translation by $f^{-1}(p_+) \in T_{p_+}S^n$. 
\item $Y^{-1}$ is translation by $(f \circ Q^{-1})(p_-) \in T_{p_-}S^n$.
\item $Y^{-1} \circ f \circ Q^{-1}$ fixes both $p_+$ and $p_-$, so it
is an element of $\Scale(p_+)\SO_n(p_+) = \Scale(p_+)\times \SO_n(p_+)$,
and so there is a unique solution to the equation $S \circ M = Y^{-1} \circ f \circ Q^{-1}$.
\end{itemize}
Given $(\pi,f_1,\cdots,f_j) \in \coba_n(j)$, let $f_i = Y_i \circ S_i \circ M_i \circ Q_i$ 
be as above. Let $Q_i(t)$ for $t \in [0,1]$ be the straight line in $\Trans({p_-})$ from 
$Q_i$ to the identity. 
Notice that $Y_i \circ S_i \circ M_i \circ Q_i(t)(p_-) = f_i(p_-)$ for all
$i$ and $t$. For all $t \in [0,1]$ there exists $r(t) \in (0,1]$ such that if  
$R_t \in \Scale(p_+)$ denotes scaling by $r(t)$ in $T_{p_-}S^n$, then
$(\pi,Y_1 \circ S_1 \circ R_t \circ M_1 \circ Q_1(t), \cdots, 
Y_j \circ S_j \circ R_t \circ M_j \circ Q_j(t)) \in \coba_n(j)$.
For each $t \in [0,1]$ choose $r(t)$ to be the maximal possible such element. 
This makes $r(t)$ a continuous function of $t$.  Observe
$(\pi,Y_1 \circ S_1 \circ R_1 \circ M_1 \circ Q_1(1), \cdots, 
Y_j \circ S_j \circ R_1 \circ M_j \circ Q_j(1)) = 
(\pi,Y_1 \circ S_1 \circ R_1 \circ M_1, \cdots, Y_j \circ S_j \circ R_1 \circ M_j) \in \frdi_n$.
\end{proof}
\end{prop}

{
\psfrag{sn}[tl][tl][1][0]{$S^n$}
\psfrag{tsn}[tl][tl][1][0]{$T_{p_-}S^n$}
\psfrag{dn}[tl][tl][1][0]{$D^n$}
\psfrag{dpn}[tl][tl][1][0]{}
\psfrag{p+}[tl][tl][1][0]{$p_+$}
\psfrag{p-}[tl][tl][1][0]{$p_-$}
\psfrag{sn1}[tl][tl][1][0]{$S^{n-1}$}
$$\includegraphics[width=10cm]{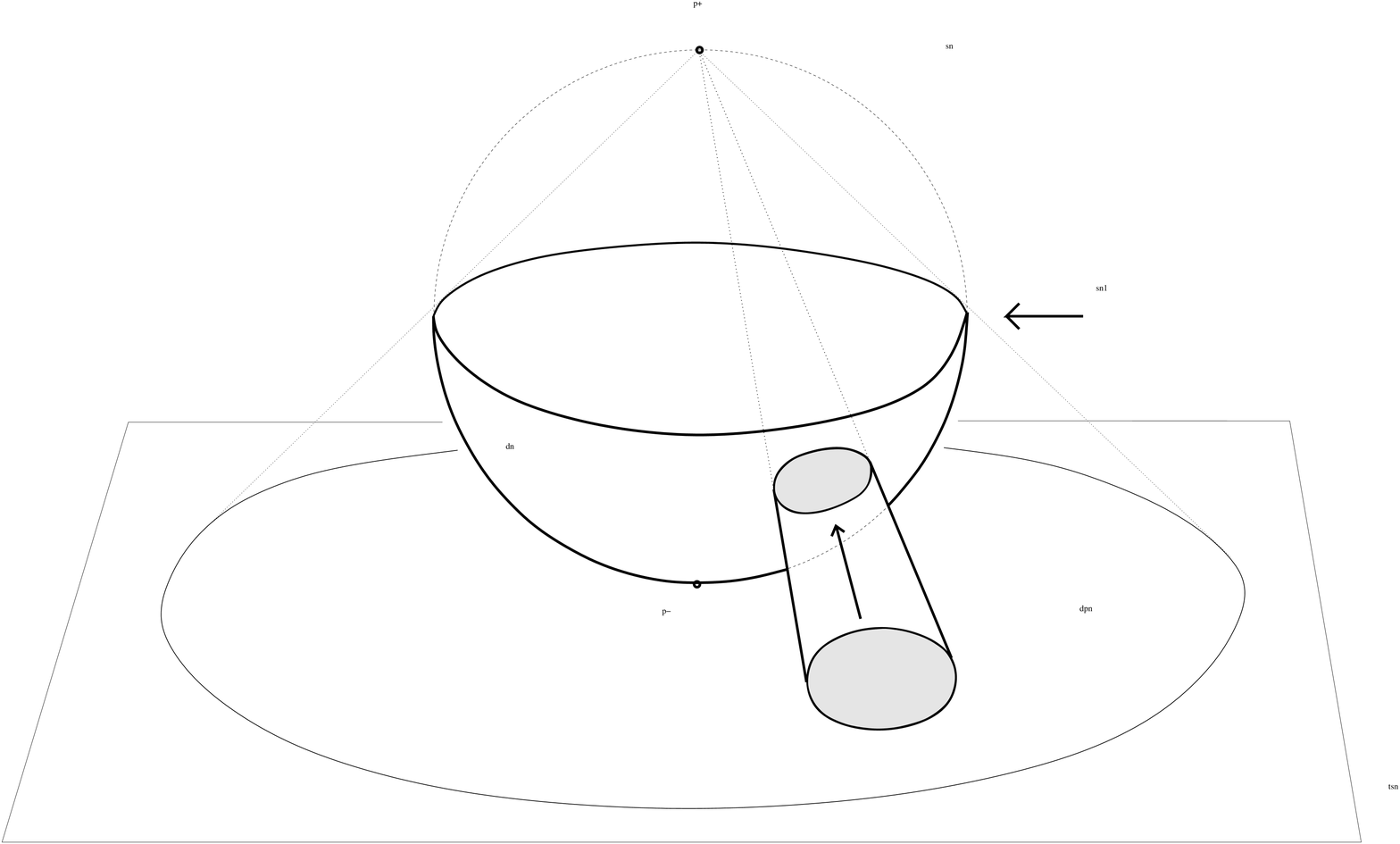}$$
\centerline{The inclusion $\frdi_n \to \coba_n$}
}

One of the main roles of the operad of framed discs $\frdi_n$
is that it detects when a space $X$ has the homotopy-type of an $n$-fold
loop space $\Omega^nY$ where $Y$ is an $\SO_n$-space \cite{PaNa}.  By the
work of Fiedorowicz \cite{Fied}, the above proposition implies that
$\coba_n$ also detects $n$-fold loop spaces over $\SO_n$-spaces.

\begin{thm}\label{mainthm} For $n \geq 2$, $\coba_n$ is a cyclic operad.
\begin{proof}
We define a right $\Sigma^+_j$-action on $\coba_n(j)$ which extends
the $\Sigma_j$ action.  Here $\Sigma^+_j$ is the permutation group of
$\{0,1,\cdots,j\}$ and $\Sigma_j$ is the stabiliser of $0$.
Given $\sigma \in \Sigma^+_j$ and $f \in \coba_n(j)$,
let $f.\sigma = \pi f^{-1}_{\sigma(0)} f \sigma$. Juxtaposition 
is interpreted as composition of functions, where $f \in \coba_n(j)$
is considered to be a an element of $\CD(S^n)^{\{0,1,\cdots,j\}}$.
More explicitly, $\sigma.(\pi,f_1,\cdots,f_j)=(\pi,\pi f^{-1}_{\sigma(0)} f_{\sigma(1)},
\cdots, \pi f^{-1}_{\sigma(0)} f_{\sigma(j)})$ where we interpret $f_0$ as $\pi$,
if it appears.
This is a right action of $\Sigma^+_j$, which can be verified by a quick
calculation:
\begin{align*} (f.\sigma).\epsilon 
 & = \left(\pi f^{-1}_{\sigma(0)} f \sigma \right).\epsilon
   = \pi f^{-1}_{\sigma\epsilon(0)} f_{\sigma(0)} \pi
	 \left( \pi f^{-1}_{\sigma(0)} 
	f \sigma \epsilon \right) \\
 & = \pi f^{-1}_{\sigma\epsilon(0)} f \sigma \epsilon
   = f.(\sigma\epsilon)
\end{align*}
The three axioms of a cyclic operad are verified below.

Axiom (i) 
$$ \xymatrix{ \coba_n(1) \ar[rr]^{\tau_1} && \coba_n(1) \\
 & 1 \ar[ul]^-{\eta} \ar[ur]_-{\eta} &
}$$
states that the transposition $(01)=\tau_1$ fixes the identity
element of $\coba_n(1)$, which can be checked: $(\pi, I).\tau_1 = 
(\pi I^{-1} I, \pi I^{-1} \pi) = (\pi,I)$.

Axiom (ii) states that the diagram below commutes:
$$ \xymatrix{
 \coba_n(i) \times \coba_n(j) \ar[r]^-{\circ_1} \ar[d]^{\tau_i \times \tau_j} & \coba_n(i+j-1) \ar[dd]^{\tau_{i+j-1}} \\
\coba_n(i) \times \coba_n(j) \ar[d]^{s} & \\
\coba_n(j) \times \coba_n(i) \ar[r]^{\circ_i} &  \coba_n(i+j-1) 
}$$
which we check by going around the diagram two ways, first the clockwise direction.
\begin{align*}
(f \circ_1 g).\tau_{i+j-1} & = (\pi,f_1 g_1,\cdots,f_1 g_j,f_2,\cdots,f_i).\tau_{i+j-1} \\
 			   & = \pi g_1^{-1} f_1^{-1} \left( f_1 g_1, f_1 g_2, \cdots, f_1 g_j, f_2, \cdots, f_i,
				 \pi \right) \\
			   & = (\pi,\pi g_1^{-1} g_2, \cdots, \pi g_1^{-1} g_j, \pi g_1^{-1} f_1^{-1} f_2,
			   \cdots, \pi g_1^{-1} f_1^{-1} f_i, \pi g_1^{-1}f_1^{-1} \pi )
\intertext{and the counter-clockwise direction}
\circ_j\left(s\left(f.\tau_i,g.\tau_j)\right)\right) & = (g.\tau_j) \circ_j (f.\tau_i) \\
	& = (\pi, \pi g_1^{-1}g_2, \cdots,\pi g_1^{-1}g_j,\pi g_1^{-1} \pi) \circ_j
	    (\pi, \pi f_1^{-1}f_2, \cdots,\pi f_1^{-1}f_i,\pi f_1^{-1} \pi) \\
	& = (\pi, \pi g_1^{-1}g_2, \cdots,\pi g_1^{-1}g_j,\pi g_1^{-1} \pi^2 f_1^{-1}f_2, \cdots,
	    \pi g_1^{-1} \pi^2 f_1^{-1}f_i, \pi g_1^{-1} \pi^2 f_1^{-1}\pi ) \\
	& = (\pi, \pi g_1^{-1}g_2, \cdots,\pi g_1^{-1}g_j,\pi g_1^{-1} f_1^{-1}f_2, \cdots,
	     \pi g_1^{-1} f_1^{-1}f_i, \pi g_1^{-1} f_1^{-1} \pi )
\end{align*}

Axiom (iii) states that the diagram below commutes:
$$ \xymatrix{
\coba_n(i) \times \coba_n(j) \ar[r]^{\circ_i} \ar[d]^{\tau_i \times 1} &  \coba_n(i+j-1) \ar[d]^{\tau_{i+j-1}} \\ 
\coba_n(i) \times \coba_n(j) \ar[r]^{\circ_{i-1}} & \coba_n(i+j-1)
}$$
which we check by going around the diagram both ways, first clockwise.
\begin{align*}
(f \circ_k g).\tau_{i+j-1} & = (\pi, f_1,\cdots,f_{k-1},f_kg_1,\cdots,f_kg_j,f_{k+1},\cdots,f_i).\tau_{i+j-1} \\
	                   & = (\pi, \pi f_1^{-1} f_2, \cdots, \pi f_1^{-1} f_{k-1}, 
	\pi f_1^{-1}f_kg_1,\cdots,\pi f_1^{-1}f_kg_j,\pi f_1^{-1}f_{k+1},\cdots,\pi f_1^{-1}f_i,
	\pi f_1^{-1} \pi)
\intertext{and counter-clockwise}
f.\tau_i \circ_{k-1} g & = ( \pi, \pi f_1^{-1} f_2, \cdots, \pi f_1^{-1} f_i, \pi f_1^{-1} \pi)\circ_{k-1} (g_1,\cdots,g_j) \\
		       & = ( \pi, \pi f_1^{-1} f_2, \cdots, \pi f_1^{-1} f_{k-1}, 
	\pi f_1^{-1}f_kg_1,\cdots,\pi f_1^{-1}f_kg_j,\pi f_1^{-1}f_{k+1},\cdots,\pi f_1^{-1}f_i,
	\pi f_1^{-1} \pi)
\end{align*}
\end{proof}
\end{thm}

The operad of little $n$-cubes (or unframed little $n$-discs) is known
not to be cyclic provided $n \geq 2$ is even, since its homology is not cyclic. A proof
appears in Proposition 3.18 of \cite{Getz}.

$\frdi_1$ and $\frdi_2$ were already known to be cyclic before the publication
of this article \cite{Getz, MSS}. The proof above can be adapted to give a coherent proof 
that $\coba_n$ is cyclic for all $n \geq 1$, provided one defines $\CD(S^1)$ to be the group of orientation-preserving isometries of hyperbolic $2$-space, restricted to the circle at 
infinity.  One defines $\pi : S^1 \to S^1$ to be rotation by $180$ degrees about the centre of the circle, and the proof proceeds as above.

\providecommand{\bysame}{\leavevmode\hbox to3em{\hrulefill}\thinspace}

\Addresses

\end{document}